\documentclass[11pt, reqno]{amsart}

\usepackage{hyperref}
\usepackage[capitalize]{cleveref}
\usepackage{multicol}
\newcommand{\pref}[1]{(\ref{#1})}
\renewcommand{\eqref}[1]{\pref{#1}}
\newcommand{\fullcref}[2]{\cref{#1}\pref{#1-#2}}
\newcommand{\fullcsee}[2]{(see \fullcref{#1}{#2})}

\DeclareMathOperator{\Cay}{Cay}

\newcommand{\normal}{\triangleleft}
\newcommand{\quot}{\overline}

\newcommand{\voltage}{\mathbb{V}}

\newcommand{\cyclic}{\mathcal{C}}
\newcommand{\ZZ}{\mathbb{Z}}

\newtheorem{lem}{Lemma}
\crefformat{lem}{Lemma~#2#1#3}
\Crefname{lem}{Lemma}{Lemmas}
\crefname{lem}{lemma}{lemmas}

\newtheorem{prop}[lem]{Proposition}
\crefformat{prop}{Proposition~#2#1#3}
\Crefname{prop}{Proposition}{Propositions}
\crefname{prop}{proposition}{propositions}

\newtheorem{thm}[lem]{Theorem}
\crefformat{thm}{Theorem~#2#1#3}
\Crefname{thm}{Theorem}{Theorems}
\crefname{thm}{theorem}{theorems}

\newtheorem{FGL}[lem]{Lemma}
\crefformat{FGL}{the Factor Group Lemma~#2\textup(#1\textup)#3}

\theoremstyle{definition}
\newtheorem{defn}[lem]{Definition}
\newtheorem{notation}[lem]{Notation}

\theoremstyle{remark}
\newtheorem{rem}[lem]{Remark}

\numberwithin{lem}{section}

 \newcounter{case}
\newenvironment{case}[1][\unskip]{\refstepcounter{case}
 \medskip \noindent \textbf{Case \thecase\ #1.}\em\ }{\unskip\upshape}
\crefname{case}{case}{cases}
\crefformat{case}{Case~#2#1#3}

 \newcounter{subcase}
\newenvironment{subcase}[1][\unskip]{\refstepcounter{subcase}
 \medskip \noindent \textbf{Subcase \thesubcase\ #1.}\em\ }{\unskip\upshape}
\crefname{subcase}{subcase}{subcases}
\crefformat{subcase}{Subcase~#2#1#3}
\numberwithin{subcase}{case}

\begin{document}

\title[On hamiltonian cycles in Cayley graphs of order $pqrs$]{On hamiltonian cycles in Cayley graphs of order \larger$p q r \mkern-1.8mu s$}

\author{Dave Witte Morris}
\address{Department of Mathematics and Computer Science,
University of Lethbridge,
4401 University Drive,
Lethbridge, Alberta, T1K~3M4, Canada}
\email{dave.morris@uleth.ca}

\begin{abstract}
Let $G$ be a finite group. We show that if $|G| = pqrs$, where $p$, $q$, $r$, and~$s$ are distinct odd primes, then every connected Cayley graph on~$G$ has a hamiltonian cycle.
\end{abstract}

\date{29 July 2021} 

{\mathversion{bold}
\maketitle
}

\section{Introduction}

\begin{defn}[cf.\ {\cite[p.~34]{GodsilRoyle}}]
If $A$ is a subset of a finite group~$G$, then the corresponding \emph{Cayley graph} $\Cay(G; A)$ is the undirected graph whose vertices are the elements of~$G$, and such that vertices $g$ and~$h$ are adjacent if and only if $g^{-1} h \in A \cup A^{-1}$, where $A^{-1} = \{\, a^{-1} \mid a \in A\,\}$.
\end{defn}

It is easy to see (and well known) that $\Cay(G; A)$ is connected if and only if $A$ is a generating set of~$G$. Several papers show that all connected Cayley graphs of certain orders are hamiltonian:

\begin{thm}[see {\cite{M2Slovenian-LowOrder,Maghsoudi-6pq,MorrisWilK} and references therein}] \label{Ham<100}
 Let $A$ be a generating set of a finite group~$G$. If\/ $|G|$ has any of the following forms\/ {\textup(}where $p$, $q$, and~$r$ are distinct primes, and $k$ is a positive integer\/{\textup)}, then\/ $\Cay(G; A)$ has a hamiltonian cycle:
 \begin{multicols}{3}

 \begin{enumerate}

 \item \label{Ham<100-kp}
 $k p$, where $k \le 47$,
 
 \item \label{Ham<100-kpq}
 $k p q$, where $k \le 7$,
 
 \item \label{Ham<100-pqr}
 $p q r$,
 
 \item \label{Ham<100-kp2}
 $k p^2$, where $\le 4$, 
 
 \item \label{Ham<100-kp3}
 $k p^3$, where $k \le 2$.
 
 \item \label{Ham<100-pk}
 $p^k$.
 
 \end{enumerate}
 \end{multicols}
\end{thm}

The purpose of this note is to add $pqrs$ to the list when it is odd:

\begin{thm} \label{pqrsodd}
If $p$, $q$, $r$, and~$s$ are distinct odd primes, then every connected Cayley graph of order $p q r s$ has a hamiltonian cycle.
\end{thm}

\begin{rem}
To remove \cref{pqrsodd}'s assumption that the primes are odd, it would suffice to show that every connected Cayley graph of order $2pqr$ is hamiltonian (where $p$, $q$, and~$r$ are distinct odd primes). An important first step in this direction was taken by F.\,Maghsoudi \cite{Maghsoudi-6pq}, who handled the case where one of the primes is~$3$.

Unfortunately, current methods do not seem to be sufficient to remove the assumption that $p$, $q$, $r$, and~$s$ are distinct. For example, it is not known that all Cayley graphs of order $9p^2$ or $3p^3$ are hamiltonian (cf.\ \cite{M2Slovenian-LowOrder}).
\end{rem}

\section{Preliminaries}

We use the following fairly standard notation.

\begin{notation}
Let $A$ be a subset of a finite group~$G$.
	\begin{enumerate}
	\item $e$ is the identity element of~$G$,
	\item $|g|$ is the order of an element~$g$ of~$G$.
	\item $G' = [G, G]$ is the \emph{commutator subgroup} of~$G$.
	\item $Z(G) = \{\, z \in G \mid \text{$zg = gz$ for all $g \in G$} \,\}$ is the \emph{centre} of~$G$.
	\item A sequence $(a_i)_{i=1}^m = (a_1,a_2,\ldots,a_m)$ of elements of $A \cup A^{-1}$ represents the walk in $\Cay(G; A)$ that visits the vertices
		\[ e, a_1, a_1 a_2, \ldots, a_1 a_2 \cdots a_m . \]
	Also, we use $a^k$ and $a^{-k}$ to represent the sequences $(a,a,\ldots,a)$ and $(a^{-1}, a^{-1}, \ldots,a^{-1})$ of length~$k$.
	\end{enumerate}
\end{notation}

Most cases of \cref{pqrsodd} are easy consequences of three known results that are collected in the following \lcnamecref{hamG'}.

\begin{thm}[{\cite[p.~257]{ChenQuimpo}, Durnberger \cite{Durnberger-p}, Morris \cite{Morris-pq}}] \label{hamG'}
Let $A$ be a generating set of a nontrivial finite group~$G$, and assume that\/ $|G|$ is odd. If either
	\begin{enumerate}
	\setcounter{enumi}{-1}
	\item \label{hamG'-1}
	$|G'| = 1$,
	or
	\item \label{hamG'-p}
	$|G'|$ is prime,
	or
	\item \label{hamG'-pq}
	$|G'|$ is the product of\/ $2$ distinct primes,
	\end{enumerate}
then $\Cay(G; A)$ has a hamiltonian cycle.
\end{thm}

\begin{rem}
The assumption that $G$ has odd order is not necessary in parts~\pref{hamG'-1} and~\pref{hamG'-p} of \cref{hamG'}.
Steps toward removing this assumption from part~\pref{hamG'-pq} were taken in \cite{GhaderpourMorris} and~\cite{Morris-2p}.
\end{rem}

For ease of reference, we record a few other known facts:

\begin{FGL}[``Factor Group Lemma'' {\cite[\S2.2]{WitteGallian-survey}}] \label{FGL}
Suppose
 \begin{itemize}
 \item $A$ is a subset of the group~$G$,
 \item $N$ is a cyclic, normal subgroup of~$G$,
 \item $C = (a_i)_{i=1}^m$ is a hamiltonian cycle in $\Cay(G/N;A)$,
 and
 \item the voltage $\voltage(C) = \prod_{i=1}^m a_i$ generates~$N$.
 \end{itemize}
 Then $(a_1,a_2,\ldots,a_m)^{|N|}$ is a hamiltonian cycle in $\Cay(G;A)$.
 \end{FGL}

\begin{lem}[{}{\cite[Lem.~2.27]{M2Slovenian-LowOrder}}] \label{NormalEasy}
 Let $A$ generate a finite group~$G$ and let $b \in A$, such that the subgroup generated by~$b$ is normal in~$G$. If\/ $\Cay \bigl( G/\langle b \rangle ; A \bigr)$ has a
hamiltonian cycle,
 and $\langle b \rangle \cap Z(G) = \{e\}$, then $\Cay(G;A)$ has a hamiltonian cycle.
 \end{lem}
 
 \begin{lem}[cf.\ {\cite[Thm.~9.4.3, p.~146]{Hall-TheoryOfGroups}} and {\cite[Lem.~2.16]{Maghsoudi-6pq}}] \label{squarefree}
Assume $G$ is a finite group of square-free order. Then
	\begin{enumerate}
	\item \label{squarefree-cyclic}
	$G'$ and $G/G'$ are cyclic,
	\item \label{squarefree-Z(G)}
	$G' \cap Z(G) = \{e\}$,
	and
	\item \label{squarefree-a}
	if $a \in G$, such that $\langle a, G' \rangle = G$, then 
		\begin{enumerate}
		\item \label{squarefree-a-commute}
		$a$ does not commute with any nontrivial element of~$G'$,
		and
		\item \label{squarefree-a-order}
		$|a| = |G/G'|$.
		\end{enumerate}
	\end{enumerate}
\end{lem}

 \begin{prop}[Maghsoudi {\cite[Prop.~3.1]{Maghsoudi-6pq}}] \label{4gen}
 Assume $G$ is a finite group, such that $|G|$ is a product of four distinct primes. If $A$ is an irredundant generating set of~$G$, such that $|A| \ge 4$, then $\Cay(G; A)$ has a hamiltonian cycle.
 \end{prop}

\section{Proof of \texorpdfstring{\cref{pqrsodd}}{the main theorem}}

Let $A$ be an irredundant generating set of a finite group~$G$, such that $|G| = p q r s$. We wish to show that $\Cay(G;A)$ has a hamiltonian cycle. 
Although some of the details are new, the line of argument here is quite standard. (See the proofs in~\cite{Maghsoudi-6pq}, for example.)

Let
	\[ \quot G = G/G' . \]
Note that $G$ is solvable by \fullcref{squarefree}{cyclic}, so $G' \neq G$.   Therefore $|G'|$ has at most~$3$ prime factors.  We may assume it has exactly~$3$, for otherwise \cref{hamG'} applies. Hence, we may assume $|G'| = qrs$ (so $|\quot G| = p$). Since $G'$ is cyclic \fullcsee{squarefree}{cyclic}, we have
	\[ G' = \cyclic_q \times \cyclic_r \times \cyclic_s , \]
where $\cyclic_n$ denotes a (multiplicative) cyclic group of order~$n$.
Then each element $\gamma$ of~$G'$ can be written uniquely in the form
	\[ \text{$\gamma = \gamma_q \, \gamma_r \, \gamma_s$,
	\ where $\gamma_q \in \cyclic_q$, $\gamma_r \in \cyclic_r$, and $\gamma_s \in \cyclic_s$}
	. \]

Let $a \in A$, such that $\langle \quot a \rangle = \quot G$. Then $|a| = p$ by \fullcref{squarefree}{a-order}. Let $b$ be another element of~$A$, so we may write $b = a^k \gamma$, where $\gamma \in G'$ (and $0 \le k < p$).

\begin{case}
Assume $A \cap G' \neq \emptyset$.
\end{case}
Then we may assume $b \in G'$. Since $G'$ is cyclic \fullcsee{squarefree}{cyclic}, and it is a standard exercise in undergraduate abstract algebra to show that every subgroup of a cyclic normal subgroup is normal, we know that $\langle b \rangle \normal G$. Let $\widehat G = G/ \langle b \rangle$. Since $|G| = pqrs$ has only~4 prime factors, we know that $|\widehat G|$ has at most~3 prime factors, so we see from part~\pref{Ham<100-kp}, \pref{Ham<100-kpq}, or~\pref{Ham<100-pqr} of \cref{Ham<100} that $\Cay(\widehat G; A)$ has a hamiltonian cycle. Since 
	\[ \langle b \rangle \cap Z(G) \subseteq G'  \cap Z(G) = \{e\} \]
\fullcsee{squarefree}{Z(G)},
we can now conclude from \cref{NormalEasy} that $\Cay(G; S)$ has a hamiltonian cycle, as desired.

\begin{case}
Assume $|A| = 2$ \textup(and $A \cap G' = \emptyset$\textup).
\end{case}
This means that $A = \{a,b\}$, so $\langle a,b \rangle = G$. This implies $\langle \gamma \rangle = G'$, so $\gamma_q$, $\gamma_r$, and~$\gamma_s$ are nontrivial. Also, we have $k \ge 1$, since $A \cap G' = \emptyset$. Let 
	\[ C = ( \quot b, \quot a^{-(k-1)}, \quot b, \quot a^{p-k-1}), \]
so $C$ is a (well known) hamiltonian cycle in $\Cay(\quot G; \quot A)$. The voltage of~$C$ is
	\[ \voltage(C) = b a^{-(k-1)} b a^{p-k-1}
		= a^k \gamma \cdot a^{-(k-1)} \cdot a^k \gamma \cdot a^{p-k-1}
		= a^k \cdot \gamma_q \gamma_r \gamma_s \cdot a \cdot \gamma_q \gamma_r \gamma_s \cdot a^{-k-1}
		. \]
Since $|G|$ is odd, we know that $a$ does not invert any nontrivial element of~$G'$, so the two occurrences of~$\gamma_q$ in this product do not cancel each other, and similarly for the occurrences of $\gamma_r$ and~$\gamma_s$. Therefore, this voltage projects nontrivially to $\cyclic_q$, $\cyclic_r$, and~$\cyclic_s$, so it generates~$G'$. Hence, \cref{FGL} provides a hamiltonian cycle in $\Cay(G;A)$.

\begin{case}
Assume $|A| \ge 3$ \textup(and $A \cap G' = \emptyset$\textup).
\end{case}
Let $c$ be a third element of~$A$. We may assume $|A| < 4$ (for otherwise \cref{4gen} applies), so $A = \{a,b,c\}$. Write $c = a^\ell \gamma'$ (with $\gamma' \in G'$). 
Since $\langle a,b,c \rangle = G$, we must have $\langle \gamma, \gamma' \rangle = G'$, but the fact that $A$ is irredundant implies $\langle \gamma \rangle \neq G'$ and $\langle \gamma' \rangle \neq G'$.
If $|\gamma| = q$, then we must have $|\gamma'| = rs$. However, this implies $\langle b, c \rangle = G$, which contradicts the fact that the generating set~$A$ is irredundant.
So $|\gamma|$ and $|\gamma'|$ must each have precisely two prime factors. Therefore, we may assume without loss of generality that $|\gamma| = qr$ and $|\gamma'| = rs$. Then
	\begin{align} \label{gammaqrs}
	\text{$\gamma = \gamma_q \, \gamma_r$ and $\gamma' = \gamma'_r \, \gamma'_s$}
	, \end{align}
and each of  $\gamma_q$, $\gamma_r$, $\gamma'_r$, and $\gamma'_s$ is nontrivial.

We may assume $k, \ell \le (p-1)/2$ (by replacing some generators with their inverses, if necessary). We may also assume, without loss of generality, that $\ell \le k$. 
(If $k \neq \ell$, this implies $\ell \le (p-3)/2$.)
Also note that $k, \ell \ge 1$, since $A \cap G' = \emptyset$.

\begin{subcase} \label{p=3}
Assume $p = 3$.
\end{subcase}
This implies $\quot a = \quot b = \quot c$, so we have the following two hamiltonian cycles in $\Cay(\quot G; A)$:
	\[ \text{$C_1 = (\quot a, \quot b, \quot c)$ \ and \ $C_2 = (\quot a, \quot c, \quot b)$} . \]
Their voltages are
	\begin{align*}
	\voltage(C_1) 
		&= abc 
		= a \cdot a \gamma \cdot a \gamma' 
		= a^2 \, \gamma_q \, \gamma_r \, a \, \gamma'_r \, \gamma'_s  \\
\intertext{and}
	\voltage(C_2) 
		&= acb
		= a \cdot a \gamma' \cdot a \gamma
		= a^2 \, \gamma'_r \, \gamma'_s \, a  \, \gamma_q \, \gamma_r
	. \end{align*}
In each of these voltages, there is nothing that could cancel the factor~$\gamma_q$ or the factor~$\gamma'_s$. So both voltages project nontrivially to $\cyclic_q$ and~$\cyclic_s$.

Therefore, we may assume that both voltages project trivially to $\cyclic_r$, for otherwise we have a voltage that projects nontrivially to all three factors, and therefore generates~$G'$, so \cref{FGL} applies. Hence
	\[ a^2 \gamma_r a \gamma'_r = \voltage(C_1)_r = e = \voltage(C_2)_r = a^2 \gamma'_r a \gamma_r , \]
so, letting $x = \gamma'_r \, \gamma_r^{-1} $, we have $a x = x a$, which means that $x$ commutes with~$a$. However, we also know from \fullcref{squarefree}{a-commute} that $a$ does not commute with any nontrivial element of~$G'$. Therefore, we must have $x = e$, which means $\gamma'_r = \gamma_r$. Also, since $a$ does not invert any nontrivial element of~$G'$ (since $a$ has odd order), we know that $\gamma_r a \gamma_r \neq a$. Therefore
	\[ \voltage(C_1)_r 
	= a^2 \cdot \gamma_r a \gamma'_r
	= a^2 \cdot \gamma_r a \gamma_r
	\neq a^2 \cdot a
	= e
	. \]
This voltage therefore projects nontrivially to all three factors of~$G'$, so it generates~$G'$. Hence, \cref{FGL} applies.

\begin{subcase} \label{case:notdistinct}
Assume $\quot a$, $\quot b$, and~$\quot c$ are not all distinct.
\end{subcase}
We may assume $\quot a = \quot c$ (which means $\ell = 1$). We may also assume $p \ge 5$, for otherwise \cref{p=3} applies.
Therefore 
	\[ p - k \ge p - (p-1)/2 = (p + 1)/2 \ge (5 + 1)/2 = 3 , \]
so we have the following two hamiltonian cycles in $\Cay(\quot G; \quot A)$:
	\[ \text{$C_1 = ( \quot b, \quot a^{-(k-1)}, \quot b, \quot c, \quot a^{p-k-2})$
	\ and \ 
	$C_2 = ( \quot b, \quot a^{-(k-1)}, \quot b, \quot a, \quot c, \quot a^{p-k-3})$} . \]
Their voltages are
	\begin{align*}
	\voltage(C_1) 
		&= b a^{-(k-1)} b c a^{p-k-2}
		= a^k \gamma \cdot a^{-(k-1)} \cdot a^k \gamma \cdot a \gamma' \cdot a^{p-k-2}
		= a^k \, \gamma_q \gamma_r \, a \, \gamma_q \gamma_r \, a \, \gamma'_r \gamma'_s \, a^{-k-2}
	\intertext{and}
	\voltage(C_2) 
		&= b a^{-(k-1)} b a c a^{p-k-3}
		= a^k \gamma \cdot a^{-(k-1)} \cdot a^k \gamma \cdot a \cdot a \gamma' \cdot a^{p-k-3}
		= a^k \, \gamma_q \gamma_r \, a \, \gamma_q \gamma_r \, a^2 \, \gamma'_r \gamma'_s \, a^{-k-3}
	. \end{align*}
Since $a$ does not invert~$\gamma_q$ (because $a$ has odd order), we see that the two occurrences of~$\gamma_q$ in these voltages cannot cancel each other. Hence, both voltages project nontrivially to~$\cyclic_q$. Also, there is nothing in either voltage that could cancel the single occurrence of~$\gamma'_s$, so both voltages also project nontrivially to~$\cyclic_s$. 

Therefore, in order to apply \cref{FGL}, it suffices to show that at least one of these voltages projects nontrivially to~$\cyclic_r$. If not, then both projections are trivial, so
	\begin{align*}
	a^k \gamma_r a \gamma_r a \cdot \gamma'_r a \cdot a^{-k-3}
	&= a^k \gamma_r a \gamma_r a \gamma'_r a^{-k-2}
	\\&= \voltage(C_1)_r 
	\\&= e 
	\\&= \voltage(C_2)_r
	\\&= a^k \gamma_r a \gamma_r a^2 \gamma'_r a^{-k-3}
	\\&= a^k \gamma_r a \gamma_r a \cdot a \gamma'_r \cdot a^{-k-3}
	, \end{align*}
which implies $\gamma'_r a = a \gamma'_r$. This contradicts the fact that $a$ does not centralize any nontrivial element of~$G'$ \fullcsee{squarefree}{a-commute}.

\begin{subcase} \label{case:ell}
Assume $\ell \neq (p-3)/2$.
\end{subcase}
We may assume that the preceding case does not apply. In particular, then $\quot b \neq \quot c$, so $k \neq \ell$, so we have $\ell \le (p - 5)/2$, which implies $k + \ell \le p - 3$. Therefore, we have the following two hamiltonian cycles in $\Cay(\quot G; \quot A)$:
	\begin{align}
	\label{C1}
	C_1 &= ( \quot b, \quot a^{-(k-1)}, \quot b, \quot c, \quot a^{-(\ell-1)}, \quot c, \quot a^{p - k - \ell - 2}) \\
\intertext{and}
	C_2 &= (\quot b, \quot a^{-(k-1)}, \quot b, \quot a, \quot c, \quot a^{-(\ell-1)}, \quot c, \quot a^{p - k - \ell - 3})
	. \end{align}
Their voltages are
	\begin{align*}
	\voltage(C_1) 
	&= b a^{-(k-1)} b c a^{-(\ell-1)} c a^{p - k - \ell - 2}
	\\& = a^k \gamma \cdot a^{-(k-1)} \cdot a^k \gamma \cdot  a^\ell \gamma' \cdot  a^{-(\ell-1)} \cdot a^\ell \gamma'  \cdot  a^{p - k - \ell - 2}
	\\& = a^k \cdot \gamma_q \gamma_r a \gamma_q \gamma_r \cdot  a^\ell \cdot \gamma'_r \gamma'_s  a \gamma'_r \gamma'_s \cdot  a^{p - k - \ell - 2} \\
\intertext{and}
	\voltage(C_2) 
	&= b a^{-(k-1)} b a c a^{-(\ell-1)} c a^{p - k - \ell - 3}
	\\& = a^k \gamma \cdot a^{-(k-1)} \cdot a^k \gamma \cdot  a \cdot a^\ell \gamma' \cdot  a^{-(\ell-1)} \cdot a^\ell \gamma' \cdot  a^{p - k - \ell - 3}
	\\& = a^k \cdot \gamma_q \gamma_r a \gamma_q \gamma_r \cdot  a^{\ell+1} \cdot \gamma'_r \gamma'_s  a \gamma'_r \gamma'_s \cdot  a^{p - k - \ell - 3}
	. \end{align*}
In each product, the two occurrences of~$\gamma_q$ cannot cancel each other, and the two occurrences of $\gamma'_s$ cannot cancel each other (because $a$ does not invert any nontrivial element of~$G'$), so both voltages project nontrivially to $\cyclic_q$ and~$\cyclic_s$.

Therefore, in order to apply \cref{FGL}, it suffices to show that at least one of these voltages projects nontrivially to~$\cyclic_r$. If not, then both projections are trivial, so, much as in \cref{case:notdistinct}, we have
	\begin{align*}
	a^k \gamma_r a \gamma_r  a^\ell \cdot \gamma'_r a \gamma'_r  a \cdot a^{p - k - \ell - 3}
	& = a^k \cdot \gamma_r a \gamma_r \cdot  a^\ell \cdot \gamma'_r a \gamma'_r \cdot  a^{p - k - \ell - 2}
	\\&= \voltage(C_1)_r 
	\\&= e 
	\\&= \voltage(C_2)_r
	\\& = a^k \cdot \gamma_r a \gamma_r \cdot  a^{\ell+1} \cdot \gamma'_r a \gamma'_r \cdot  a^{p - k - \ell - 3}
	\\& = a^k \gamma_r a \gamma_r  a^\ell \cdot a \gamma'_r a \gamma'_r \cdot  a^{p - k - \ell - 3}
	, \end{align*}
so $\gamma'_r a \gamma'_r  a = a \gamma'_r a \gamma'_r$. 

Now, note that $\langle \gamma'_r \rangle \normal G$ (since every subgroup of the cyclic normal subgroup~$G'$ is normal), so there is some $\widehat\gamma \in \langle \gamma'_r \rangle$, such that $a \gamma'_r = \widehat\gamma a$. With this notation, the conclusion of the preceding paragraph tells us that 
	$\gamma'_r \widehat\gamma a^2 = \widehat\gamma a^2 \gamma'_r$. 
Since $\gamma'_r$ and $\widehat\gamma$ are in the abelian group $\langle \gamma'_r \rangle$, we know that they commute with each other, so this implies $\gamma'_r a^2 = a^2 \gamma'_r$. 

However, since $\quot a$ generates~$\quot G$, and $|\quot G| = p$ is odd, we know that $\quot a^2$ also generates~$\quot G$. Hence, we see from \fullcref{squarefree}{a-commute} that $a^2$ does not centralize any nontrivial element of~$G'$. This contradicts the conclusion of the preceding paragraph.

\begin{subcase} \label{p=7}
Assume $p = 7$.
\end{subcase}
We may assume \cref{case:notdistinct} does not apply, so $1 < \ell < k \le (p - 1)/2$. Since $(p-1)/2 = (7 - 1)/2 = 3$, we conclude that $\ell = 2$ and $k = 3$.
Thus, we have
	\[ \text{$p = 7$, \ $\quot c = \quot a^2$, \ and \ $\quot b = \quot a^3$.} \]

Now, the Cayley graph $\Cay(\quot G; \quot a, \quot b, \quot c)$ is the complete graph~$K_7$, so it has \emph{many} hamiltonian cycles.
In particular, we have the hamiltonian cycle
	\[ C = (\quot c, \quot a^2, \quot c, \quot a^{-1}, \quot b, \quot a^{-1} ) . \]
Since $G'$ is a cyclic, normal subgroup of~$G$ (and $|a| = p = 7$), there is some $\alpha \in \ZZ^+$, such that $\alpha^7 \equiv 1 \pmod{r}$ and
	\[ \text{$a \, g = g^\alpha \, a$ \ for all $g \in G'$} . \]
Also, we may write $\gamma'_r = \gamma_r^x$ for some $x \in \ZZ^+$.
With this notation, the voltage of~$C$ is
	\begin{align*}
	\voltage(C)
	&= c  a^2  c a^{-1}  b  a^{-1}  
	\\&= a^2 \gamma'  \cdot a^2  \cdot a^2 \gamma' \cdot a^{-1}  \cdot a^3 \gamma \cdot  a^{-1}
	\\&= a^2 \gamma'  a^4 \gamma' a^2 \gamma a^{-1}
	\\&= (\gamma')^{\alpha^2 + \alpha^6} \gamma^\alpha
	\\&= (\gamma'_r \gamma'_s)^{\alpha^2+ \alpha^6} (\gamma_q \, \gamma_r)^{\alpha} 
	\\&= (\gamma_r^x \gamma'_s)^{\alpha^2+ \alpha^6} (\gamma_q \, \gamma_r)^{\alpha} 
	\\&= \gamma_q^{\alpha} \, \gamma_r^{(\alpha^2 + \alpha^6)x + \alpha} \, (\gamma'_s)^{\alpha^2+ \alpha^6}
	. \end{align*}
Since $\alpha^7 \equiv 1 \pmod{r}$, it is easy to see that neither $\alpha$ nor~$\alpha^2 + \alpha^6$ is divisible by~$r$, so it is clear that $\voltage(C)$ projects nontrivially to $\ZZ_q$ and $\ZZ_s$.

Therefore, if \cref{FGL} does not apply, then this voltage projects trivially to~$\cyclic_r$, which means
	\begin{align} \label{alpha}
	 (\alpha^2 + \alpha^6)x + \alpha \equiv 0 \pmod{r} 
	 , \end{align}
We also have the hamiltonian cycle~$C_1$ of Equation~\eqref{C1}. If its voltage also does not generate~$G'$, then this voltage, too, projects trivially to~$\cyclic_r$, which means 
	\begin{align*}
	e
	&= \voltage(C_1)_r
	\\&= a^k \cdot \gamma_r a \gamma_r \cdot  a^\ell \cdot \gamma'_r a \gamma'_r  \cdot  a^{p - k - \ell - 2}
	\\&= a^3 \cdot \gamma_r a \gamma_r \cdot  a^2 \cdot \gamma_r^x a \gamma_r^x
	\\&= \gamma_r^{\alpha^3 + \alpha^4 + \alpha^6 x + \alpha^7 x}
	, \end{align*}
so 
	\[ \alpha^3 + \alpha^4 + \alpha^6 x + \alpha^7 x \equiv 0 \pmod{r} . \]
Multiplying by $\alpha^4$, and recalling that $\alpha^7 \equiv 1 \pmod{r}$, we conclude that
	\[ (1 + \alpha) + \alpha^3(1 + \alpha)x  \equiv 0 \pmod{r} . \]
Dividing this equation by $1 + \alpha$ yields $1 + \alpha^3 x = 0$ (in~$\ZZ_r$), so $x = -1/\alpha^3$.
Plugging this into Equation~\eqref{alpha} yields
	$ - (\alpha^2 + \alpha^6)/\alpha^3 + \alpha = 0 $ in~$\ZZ_r$,
so $\alpha^4 - \alpha^2 + 1 \equiv 0 \pmod{r}$. Recall that we also know $\alpha^7 - 1 \equiv 0 \pmod{r}$. We therefore have
	\begin{align*}
	0
	&= (\alpha^6 - \alpha^5 + \alpha^4 - \alpha^3 - 1) \cdot 0 \ - \ (\alpha^3 - \alpha^2) \cdot 0
	\\&\equiv (\alpha^6 - \alpha^5 + \alpha^4 - \alpha^3 - 1)(\alpha^4 - \alpha^2 + 1) - (\alpha^3 - \alpha^2)(\alpha^7 - 1) 
		&& \pmod{r}
	\\&= -1 
	. \end{align*}
This is a contradiction.

\begin{subcase}
Assume that the preceding cases do not apply.
\end{subcase}
Since \cref{case:ell} does not apply, we must have $\ell = (p-3)/2$. Since $\ell < k \le (p-1)/2$, this implies 
	\[ k = \frac{p-1}{2} = 1 + \frac{p-3}{2} = 1 + \ell . \]
Therefore 
	\[ \quot a \, \quot c =  \quot a \, \quot a^\ell = \quot a^{1 + \ell} = \quot a^k = \quot b . \]

Putting $c$ into the role of~$a$ will usually give us new values for $k$ and~$\ell$. Indeed, if \cref{case:ell} does not apply after this change, then either
	\[ \text{$\quot a^{\pm 1} = c^{(p-3)/2}$ \ and \ $\quot b^{\pm 1} = c^{(p-1)/2} = \quot a^{\pm 1} \, \quot c$} \]
or 
	\[ \text{$\quot b^{\pm 1} = c^{(p-3)/2}$ \ and \ $\quot a^{\pm 1} = c^{(p-1)/2} = \quot b^{\pm 1} \, \quot c$} .\]
However, we have
	\[ \quot a^{-1} \, \quot c = \quot a^{\ell - 1} = \quot a^{(p - 5)/2} \notin \{\quot a^{\pm(p-1)/2}\} = \{\quot b^{\pm 1}\} \]
	and
	\[ \quot b \, \quot c = \quot a^{k + 1} = \quot a^{(p + 1)/2} \notin \{\quot a^{\pm1}\} .\]
Therefore, $\quot c^{(p-3)/2}$ must be either $\quot a$ or~$\quot b^{-1}$. Noting that (since $\quot a \, \quot c = \quot b$) we have $\quot b = \quot a \, \quot c$ and $\quot a^{-1} = \quot b^{-1} \, \quot c$, this implies that the only possibilities are that either
	\begin{align} \label{ab}
	\text{$\quot a = \quot c^{(p-3)/2}$ \ and \ $\quot b = \quot a \, \quot c = \quot{c}^{(p-1)/2}$} 
	\end{align}
or 
	\begin{align} \label{binvainv}
	\text{$\quot b^{-1} = \quot c^{(p-3)/2}$ \ and \ $\quot a^{-1} = \quot b^{-1} \, \quot c = \quot{c}^{(p-1)/2} $}
	. \end{align}

If \eqref{ab} holds, then 
	\[ \quot a = \quot c^{(p-3)/2} = \bigl( \quot a^{(p-3)/2} \bigr)^{(p-3)/2} = \quot a^{(p-3)^2/4}, \]
so $(p-3)^2/4 \equiv 1 \pmod{p}$, so $9 \equiv 4 \pmod{p}$, which implies $p = 5$. 
However, since \cref{case:notdistinct} does not apply, we know that $\quot c \neq \quot a$, so $\ell > 1$. This means $(p - 3)/2 \ge 2$, so $p \ge 7$. This is a contradiction.

So \eqref{binvainv} must hold. Then $\quot b^{-1} = \quot c^{(p-3)/2}$, so we must have 
	\[ (p-3)^2/4 \equiv -(p-1)/2 \pmod{p} . \]
so $9 \equiv 2 \pmod{p}$, which implies $p = 7$. So \cref{p=7} applies.

\end{document}